\documentclass{article}

\usepackage{latexsym, amsmath, amssymb, amsthm}

\newtheorem{theorem}{Theorem}
\newtheorem{proposition}[theorem]{Proposition}
\newtheorem{lemma}[theorem]{Lemma}
\newtheorem{corollary}[theorem]{Corollary}
\theoremstyle{definition}

\theoremstyle{remark}
\newtheorem*{remark}{Remark}

\newcommand{\norm}[1]{\left\| #1 \right\|}
\newcommand{\N}{\mathbb{N}}

\def\R{\mathbb{R}} 
\def\Z{\mathbb{Z}} 
\def\F{\mathcal{F}}

\def\vol{\mathrm{Leb}}

\def\lap{\ell} 
\def\t{\mathbf{t}}
\def\cC{\mathcal{C}}%
\def\C{\mathbf{C}}
\def\D{\mathcal{D}}
\def\sym{\mathcal{A}}
\def\a{\mathbf{a}}
\def\b{\mathbf{b}}
\def\c{\mathbf{c}}
\def\x{\mathbf{x}}
\def\y{\mathbf{y}}
\def\P{\mathcal{P}}

\def\e{\mathbf{e}}
\def\ps{S^1\times\R}
\def\ann{D}
\def\ac{{\dagger}}
\DeclareMathOperator{\supp}{supp}
\DeclareMathOperator{\Jac}{\mathbf{Jac}}
\DeclareMathOperator{\Ker}{\mathrm{Ker}}
\newcommand{\boB}{\mathcal{B}}
\newcommand{\boC}{\mathcal{C}}
\newcommand{\tq}{\ |\ }
\DeclareMathOperator{\inte}{Int}
\newcommand{\ic}{\mathbf{i}}
\newcommand{\ru}{\rho_0}
\newcommand{\rv}{\rho_1}
\DeclareMathOperator{\Int}{Int}

\title{Smoothness of solenoidal attractors}
\author{
Artur Avila\thanks{
\emph{Address:} Laboratoire de Probabilit\'es et Mod\`eles
al\'eatoires,
Universit\'e Pierre et Marie Curie--Bo\^\i te courrier 188,
75252--Paris Cedex 05, France; \emph{Email:} artur@ccr.jussieu.fr},
S\'ebastien Gou\"ezel\thanks{ \emph{Address:} D\'epartement de
Math\'ematiques et Applications, Ecole Normale Sup\'erieure, 45 rue
d'Ulm, Paris, France; \emph{Email:} Sebastien.Gouezel@ens.fr} and
Masato Tsujii\thanks{ \emph{Address:} Department of Mathematics, Hokkaido University, Kita 10 Nishi 8, Kita-ku, Sapporo, Japan; \emph{Email:} tsujii@math.sci.hokudai.ac.jp }
}

\date{February 24, 2005}

\begin{document}

\maketitle

\begin{abstract}
We consider dynamical systems generated by skew products of affine
contractions on the real line over angle-multiplying maps on the
circle $S^1$:
\[
T:S^{1}\times \R\to S^{1}\times \R,\qquad T(x,y)=(\lap x, \lambda
y+f(x))
\]
where $\lap\ge 2$, $0<\lambda<1$ and $f$ is a $C^{r}$ function on
$S^{1}$. We show that, if $\lambda^{1+2s}\lap>1$ for some $0\leq
s< r-2$, the density of the SBR measure for $T$ is contained
in the Sobolev space $W^{s}(S^{1}\times \R)$ for almost all ($C^r$
generic, at least) $f$.
\end{abstract}

\section{Introduction}
In this paper, we study dynamical systems generated by skew
products of affine contractions on the real line over
angle-multiplying maps on the circle $S^1=\R/\Z$:
\begin{equation}\label{eqn:Tdef}
 T:S^{1}\times \R\to
S^{1}\times \R,\qquad
 T(x,y)=(\lap x, \lambda y+f(x))
\end{equation}
where $\lap\ge 2$ is an integer, $0<\lambda<1$ is a real number
and $f$ is a $C^{r}$ function on~$S^{1}$ (for some integer $r \geq 3$).
It admits a forward
invariant closed subset $A$ such that $\omega(\x)=A$ for Lebesgue
almost every point $\x\in \ps$. Further, there exists an ergodic
invariant probability measure $\mu$ on $A$ for which Lebesgue
almost every point on $\ps$ is generic. The measure $\mu$ is
called the {\em SBR measure} for $T$. If $T$ is locally area
contracting, {\it i.e.}, $\det DT\equiv \lambda \lap<1$,  the
subset $A$ is a Lebesgue null subset and hence the SBR measure is
totally singular with respect to the Lebesgue measure. In
\cite{T}, the third named author studied the case where $T$ is
locally area expanding, {\it i.e.}, $\lambda \lap >1$, and proved
that the SBR measure is absolutely continuous with respect to the
Lebesgue measure for $C^r$ generic $f$.

In the present paper, we study the smoothness of the density of
the SBR measure in more detail, and the mixing properties of $T$.

\begin{theorem}\label{th:main1}
If $\ell$ and $\lambda$ satisfy $\lambda^{1+2s} \lap>1$ for some
$0\leq s< r-2$, the density of the SBR measure $\mu$ with
respect to the Lebesgue measure is contained in the Sobolev space
$W^{s}(\ps)$ for any $f$ in an open dense subset of $C^r(S^1,\R)$.
\end{theorem}
Since the elements of $W^{s}(\ps)$ for $s>1$ are continuous up to
modification on Lebesgue null subsets from Sobolev's embedding
theorem, it follows
\begin{corollary}
If $\lambda^{3} \lap>1$ and $r\ge 4$, the attractor $A$ has
non-empty interior for  $f$ in an open dense subset of
$C^r(S^1,\R)$.
\end{corollary}
\begin{remark}
Recently, Bam\'on, Kiwi and  Rivera-Letelier announced the
following result: for an open dense subset of $C^{1+ \epsilon}$
hyperbolic endomorphisms of the annulus, $\log d+\chi>0$ implies
that the attractor has non-empty interior, where  $d$  is the
degree of the induced map in homology and  $\chi$ is the negative
Lyapunov exponent of the SBR measure. (See also \cite{B}.)
\end{remark}
\begin{remark}
When $s>1$, we also obtain that the density of the SBR measure is
essentially bounded. Together with the results of Rams in \cite{R}, it
gives examples of solenoids in higher dimensions for which the
invariant measure is equivalent to the Hausdorff measure.
\end{remark}

The Perron-Frobenius operator $P:L^1(\ps)\to L^1(\ps)$ is defined
by
  \[
  Ph(\x)=\frac{1}{\lambda\lap}\sum_{\y\in T^{-1}(\x)} h(\y),
  \]
and characterized by the  property that
\begin{equation}\label{eqn:ch}
\frac{dT_*\nu}{d\vol}=P\left(\frac{d\nu}{d\vol}\right)
\end{equation}
for any finite measure $\nu$ which is  absolutely continuous with respect to
the Lebesgue measure $\vol$ on $\ps$.

When $s>1/2$, we obtain a precise spectral description of $P$, which
strengthens considerably Theorem~\ref{th:main1}.
\begin{theorem}
\label{th:main2}
Assume that $\ell$ and $\lambda$ satisfy $\lambda^{1+2s} \lap>1$
for some $1/2< s< r-2$. Take $\gamma\in
\left ((\lambda^{1+2s}\lap)^{-1/2},1 \right )$.
For any $f$ in an open dense subset of
$C^r(S^1,\R)$,
there exists a Banach space $\boB$ contained in $W^s(\ps)$
on which the transfer operator $P$ acts continuously with
an essential spectral radius at most $\gamma$ (in particular, $P$ admits a
spectral gap, and the correlations of $T$ decay exponentially fast).
Moreover, $\boB$ can be chosen to contain
all functions in $C^{r-1}(\ps)$ supported in some given (fixed) compact
subset of $\ps$.
\end{theorem}
Since $T$ is uniformly hyperbolic, the exponential decay of
correlations was already known. The novel feature of our theorem
is that, when the contraction coefficient $\lambda$ tends to $1$,
our estimates do \emph{not} degenerate. In fact, the inequality
$\lambda<1$ is used only to ensure that a compact subset of $S^1
\times \R$ is invariant, to get an SBR measure. Hence, our method
may probably be generalized to settings with a neutral (or
slightly positive) exponent on a compact space.

Fix $\lap\ge 2$ and let $\D_{r,s}\subset (0,1)\times
C^{r}(S^{1},\R)$ be the set of pairs $(\lambda,f)$ such that the
conclusions of Theorems~\ref{th:main1} and~\ref{th:main2} hold.  Let
$\D^{\circ}_{r,s}$ be the interior of $\D_{r,s}$. The following
result shows that Theorems~\ref{th:main1} and~\ref{th:main2}
hold for ``almost all''
$T$, in a precise sense:
\begin{theorem}\label{th:main}
If  $\ell$ and $\lambda$ satisfy $\lambda^{1+2s} \lap>1$ for some
$0\leq s< r-2$, there exists a finite collection of
$C^{\infty}$ functions $\varphi_{i}:S^{1}\to \R$, $1\le i\le m$,
such that, for any $g\in C^{r}(S^{1},\R)$, the subset
\[
\left\{(t_{1},t_{2},\cdots,t_{m})\in \R^{m}
\;\left|\;\left(\lambda,\; g(x)+\sum_{i=1}^{m}
t_{i}\varphi_{i}(x)\right)\notin \D^{\circ}_{r,s}\right.\right\}
\]
is a null subset with respect to the Lebesgue measure on $\R^{m}$.
\end{theorem}

We proceed as follows. In the next section, we introduce some
definitions related to a transversality condition on the mapping
$T$, which is similar to (but slightly different from) that used
in \cite{T}. This transversality condition is proved to be a
generic one in the last section. In Section~\ref{sec:pf}, we
introduce some norms on the space of $C^r$ functions on $\ps$ and
prove a Lasota-Yorke type inequality for them, imitating the
argument in the recent paper \cite{GL} of C. Liverani and the
second named author with slight modification. Section~\ref{sec:ly}
is the core of this paper, where we prove a Lasota-Yorke inequality
involving the $W^s$ norm and the norm introduced in
Section~\ref{sec:pf}. Finally, in Section~\ref{sec:proofs_thms}, we show how
these Lasota-Yorke inequalities imply the main results of the paper.

\section{Some definitions}\label{sec:def}
 From here to the end of this paper, we fix an integer $\lap\ge 2$,
real numbers $0<\lambda<1$ and $0\le s <r-2$ satisfying
$\lambda^{1+2s}\lap>1$. We also fix a positive number $\kappa$ and
consider the mapping $T$ for a function $f$ in
  \[
  \mathcal{U}=\mathcal{U}_\kappa=
  \left\{ f\in C^r(S^1,\R)\;;\; \|f\|_{C^r}:=\max_{0\le
  k\le r} \sup_{x\in S^1}\left|\frac{d^k}{dx^k}f(x)\right|\le \kappa
  \right\}.
  \]
Fix ${\alpha_0} =\kappa/(1-\lambda)$ and let $\ann=S^1\times
[-{\alpha_0},{\alpha_0}]$. Then we have $T(\ann)\subset \ann$. Let
$\P$ be the  partition of $S^{1}$ into the intervals $\P(k)=[
(k-1)/\lap,  k/\lap)$ for $1\le k\le \ell$.
Let $\tau:S^{1}\to S^{1}$ be the map defined by $\tau(x)=\lap
\cdot x$.
Then the partition
$\P^{n}:=\bigvee_{i=0}^{n-1}\tau^{-i}(\P)$ for $n\ge 1$ consists
of the intervals
\[
\P(\a)=\bigcap_{i=0}^{n-1}\tau^{-i}\left(\P(a_{n-i})\right),
\qquad \a=(a_{i})_{i=1}^{n}\in\sym^{n}
\]
where $\sym^{n}$ denotes the space of words of length $n$ on the
set $\sym=\{1,2,\cdots,\lap\}$.
\begin{remark}Notice that $\a$ is the {\em reverse} of the
itinerary of points in $\P(\a)$.
\end{remark}
For $x\in S^{1}$ and $\a\in \sym^{n}$, there is a unique
point $y\in \P(\a)$ such that $\tau^{n}(y)=x$, which is denoted by
$\a (x)$. For $\a=(a_{i})_{i=1}^{n}\in \sym^{n}$, the image of the
segment $\P(\a)\times\{0\}\subset S^{1}\times \R$ under the
iterate $T^{n}$ is the graph of the function $S(\cdot,\a)$ defined
by
\begin{equation*}
S(x,\a):=\sum_{i=1}^{n}\lambda^{i-1}f(\tau^{n-i}(\a(x)))
=\sum_{i=1}^{n}\lambda^{i-1}f([\a]_{i}(x))
\end{equation*}
where $[\a]_{q}=(a_{i})_{i=1}^{q}$. For a word
$\a=(a_{i})_{i=1}^{\infty}\in \sym^{\infty}$ of infinite length,
we define
\begin{equation*}
S(x,\a)=\lim_{i\to\infty}S(x,[\a]_{i})=\sum_{i=1}^{\infty}\lambda^{i-1}f([\a]_{i}(x)).
\end{equation*}
For a word $\c$ of  length $m$, let $\P_*(\c)$ be the union of the
interval $\P(\c)$ and the two intervals in $\P^{m}$ adjacent to
it. The function $S(\cdot,\a)$ for a word $\a\in \sym^n$ with
$1\le n\le \infty$ may not be continuous on $\P_*(\c)$ when
$\P(\c)$ has $0\in S^1$ as its end. Nevertheless  the restriction
of $S(\cdot,\a)$ to $\P(\c)$ can be naturally extended to
$\P_*(\c)$ as a $C^{r}$ function. Indeed, letting
$\tau^{-i}_{\c,\a}:\P_*(\c)\to S^1$ be the branch of the inverse
of $\tau^i$ satisfying  $\tau^{-i}_{\c,\a}(\P(\c))\subset
\P([\a]_i)$, the extension is given by
\begin{equation} \label{eqn:Sc}
S_{\c}(\cdot,\a) : \P_*(\c)\to\R,\quad
S_{\c}(x,\a):=\sum_{i=1}^{n}\lambda^{i-1}f(\tau^{-i}_{\c,\a}(x)).
\end{equation}
For any word $\a$ of finite or infinite length, we have
\begin{equation}\label{eqn:alpha}
\sup_{x\in \P_*(\c)}\max_{0\le \nu\le r} \lap^\nu \left|\frac{d^{\nu}}{d
x^{\nu}}S_{\c}(x,\a)\right|\le {\alpha_0}.
\end{equation}

For $\a,\b\in \sym^{q}$ and $\c\in \sym^{p}$, we say that $\a$ and
$\b$ are {\em transversal} on $\c$ and write $\a\pitchfork_{\c}\b$
if
  \[
  \left|\frac{d}{d x}S_{\c}(x,\a)-\frac{d}{d
  x}S_{\c}(y,\b)\right|>2\lambda^q\ell^{-q} {\alpha_0}
  \]
at all points $x,y$ in the closure of $\P_*(\c)$. We put
  \[
  \e(q,p)=\max_{\c\in\sym^{p}}\max_{\a\in \sym^{q}}
  \#\{\b\in\sym^{q}\mid \a\not \pitchfork_{\c} \b\} \qquad \mbox{
  and}\qquad \e(q)=\lim_{p\to\infty}\e(q,p).
  \]

The main argument of the proof will be to construct norms which will
satisfy a Lasota-Yorke inequality if $\e(q)$ is not too big for
some $q$. This will readily imply the two main theorems if the norms
have sufficiently good properties. To conclude, a transversality
argument (similar to the arguments in \cite{T}) will show that, for
almost all functions $f$ (in the sense of Theorem \ref{th:main}),
$\e(q)$ is not too big for some $q$.

Henceforth, and until the end of Section \ref{sec:ly},
we fix a large integer $q$. By definition,
there exists $p_0\geq 1$ such that $\e(q,p)=\e(q)$ for $p\ge p_0$.
We also fix an integer $p\geq p_0$.

\section{Perron-Frobenius operator and the norm $\|\cdot\|^\ac_\rho$}
\label{sec:pf}

Let $C^r(\ann)$ be the set of $C^r$ functions on $\ps$ whose
supports are contained in $\ann$. In this section, we define
preliminary norms
on the space $C^r(\ann)$ and show Lasota-Yorke type inequalities
for them. For the definition of the norms, we prepare a class
$\Omega$ of $C^r$ curves on $\ps$. Let
$\gamma:\mathcal{D}(\gamma)\to \ps$ be a continuous curve on $\ps$
whose domain of definition $\mathcal{D}(\gamma)$ is a compact
interval. For $n\ge 0$, there are $\ell^n$ curves
$\tilde\gamma_i:\mathcal{D}(\gamma) \to
\ps$, $1\le i\le \ell^n$, such that $T^n\circ
\tilde{\gamma}_i=\gamma$, each of which is called a backward
image of $\gamma$ by $T^n$.  From the hyperbolic properties of $T$,
we can choose positive constants $c_i$, $1\le i\le r$, so that the
following holds: Let $\Omega$ be the set of $C^r$ curves
$\gamma:\mathcal{D}(\gamma)\to \ps$ such that
\begin{itemize}
\item the domain of definition $\mathcal{D}(\gamma)$ is a compact interval,
\item $\gamma$ is written in the form $\gamma(t)=(\pi\circ \gamma(t),t)$ and
\item $|d^i(\pi\circ \gamma)/dt^i(s)|\le c_i$ for $1\le i\le r$ and
$s\in \mathcal{D}(\gamma)$
\end{itemize}
where $\pi:\ps\to S^1$ is the projection to the first component.
Then each backward image $\tilde\gamma$ of any $\gamma\in \Omega$
by $T^n$ with $n\ge 1$ is the composition $\hat{\gamma}\circ  g$
of a curve $\hat{\gamma}\in \Omega$ and a $C^r$ diffeomorphism $
g:\mathcal{D}(\gamma)\to \mathcal{D}(\hat{\gamma})$. Further, we
can take a positive constant $c$ so that the diffeomorphism $ g$
always satisfies
  \begin{equation}\label{eqn:distg}
  \left|\frac{d^\nu}{ds^\nu}(
  g^{-1})(s)\right|<c\lambda^{n}\quad\mbox{ for  $s\in
  \mathcal{D}(\hat{\gamma})$ and $1\le \nu \le r$.}
  \end{equation}
We henceforth fix such $c$, $c_i$, $1\le i\le r$, and  $\Omega$ as
above. Moreover, the cone
  \begin{equation}
  \label{eq:defC}
  \C=\{(u,v) \tq |u|\leq \alpha_0^{-1}|v|\}
  \end{equation}
is invariant under $DT^{-1}$, whence we can take
$c_1=\alpha_0^{-1}$. Finally, increasing the constants $c_2,\dots,c_r$
if necessary, we can assume that, whenever $I$ is a segment in
$S^1\times \R$ and $J$ is a component of $T^{-q}(I)$ such that its
tangent vectors are all contained in $\C$, then $J$ is the image of an
element of $\Omega$ (recall that $q$ is fixed once and for all until
the end of Section \ref{sec:ly}).

For a function $h\in C^r(\ann)$ and an integer $0\le \rho\le r-1$,
we define
  \[
  \|h\|^\ac_{\rho}:=\max_{{\alpha}+\beta\le \rho}\;\sup_{\gamma\in
  \Omega}\;\;\sup_{ \varphi\in \cC^{{\alpha}+\beta}(\gamma)} \int
  \varphi(t) \cdot \partial_x^{\alpha} \partial_y^\beta h(\gamma(t))
  dt
  \]
where $\max_{{\alpha}+\beta\le \rho}$ denotes the maximum over
pairs $({\alpha},\beta)$ of non-negative integers such that
${\alpha}+\beta\le \rho$ and $\cC^{s}(\gamma)$ denotes the space
of $C^s$ functions $\varphi$ on $\R$ such that $\supp \varphi
\subset \Int(\mathcal{D}(\gamma))$ and $\|\varphi\|_{C^s}\le 1$. This is
a norm on $C^r(\ann)$. It satisfies
  \begin{equation}
  \label{ac_greater_L1}
  \norm{h}_{L^1} \leq C \norm{h}_0^\ac\leq C\norm{h}_\rho^\ac.
  \end{equation}
The following lemma is the main ingredient
of this section.
\begin{lemma}\label{lm:P}
There exists a constant $A_0$ such that
  \begin{align}\label{eqn:GLLY1}
  &\|P^{n} h\|^\ac_{\rho}\le  A_0 \ell^{-\rho n} \|h\|^\ac_{\rho}
  +C(n)\|h\|^\ac_{\rho-1}\quad\mbox{ for $1\le \rho\le r-1$},
  \intertext{and}\label{eqn:GLLY2} &\|P^{n} h\|^\ac_{0}\le  A_0
  \|h\|^\ac_{0}
  \end{align}
for $n\ge 0$ and $h\in C^r(\ann)$, where  $C(n)$ may depend on $n$
but not  on $h$.
\end{lemma}
\begin{proof}
Note that the iterate $T^n$ for $n\ge 0$ is locally written in the
form
\begin{equation}\label{eqn:Tloc}
T^n(x,y)=(\ell^nx, \lambda^n y+S(\ell^n x))
\end{equation}
where $S$ is a $C^r$ function whose derivatives up to order $r$ are
bounded by ${\alpha_0}$. Consider non-negative integers $\rho$,
$\alpha$, $\beta$ satisfying $1\le \rho \le r-1$ and
$\alpha+\beta=\rho$.  Differentiating both sides of
  \[
  P^n h(x,y)=\frac{1}{\lambda^{n}\lap^n}\sum_{(x',y')\in
  T^{-n}(x,y)} h(x',y')
  \]
by using (\ref{eqn:Tloc}), we see that the differential
$\partial_x^{\alpha}\partial_y^\beta P^n h(x,y)$ can be written as
the sum of
  \begin{align*}
  \Phi(x,y)&=\sum_{(x',y')\in T^{-n}(x,y)} \sum_{k=0}^{\alpha}Q_k(x)
  \frac{\partial_x^{{\alpha}-k}\partial_y^{\beta+k} h(x',y')}
  {\lambda^{(1+\beta+k)n}\ell^{(1+\alpha-k)n}} \intertext{and}
  \Psi(x,y)&=\sum_{(x',y')\in T^{-n}(x,y)}\sum_{a+b\le
  \rho-1}Q_{a,b}(x)\frac{\partial_x^a \partial_y^b h(x',y')}
  {\lambda^{(1+b)n}\ell^{(1+a)n}}
  \end{align*}
where $Q_k(\cdot)$ and $Q_{a,b}(\cdot)$ are functions of class
$C^\rho$ and $C^{a+b}$ respectively.\footnote {Strictly speaking, the
functions $Q_k$ and $Q_{a,b}$ are only defined on the $\lap^n$-fold
covering of
$S^1$, since their definition involves the choice of an inverse branch, but
we will keep the dependence on the choice of the inverse branch implicit in
the notation.} It is easy to check that the
$C^\rho$ norm of $Q_k(\cdot)$ and $C^{a+b}$ norm of $Q_{a,b}(\cdot)$
are  bounded by some constant.

For $\gamma\in \Omega$ and $\varphi\in \cC^\rho(\gamma)$, we
estimate
  \begin{equation}\label{eqn:p}
  \int\varphi(t)\partial_x^{\alpha}\partial_y^\beta P^n h(\gamma(t))
  dt =\int\varphi(t)\Phi(\gamma(t)) dt+\int\varphi(t)\Psi(\gamma(t))
  dt.
  \end{equation}
Let $\gamma_i$, $1\le i\le \ell^n$, be the backward images of the
curve $\gamma$ by $T^n$ and write them as the composition
$\hat{\gamma}_i\circ  g_i$ of $\hat{\gamma}_i\in \Omega$ and a
$C^r$ diffeomorphism $ g_i$. Then we have
  \begin{align*}
  &\int \varphi(t)  \Psi(\gamma(t)) dt = \sum_{1\le i\le
  \ell^n}\sum_{a+b\le \rho-1} \int \varphi(t)
  \frac{Q_{a,b}(\pi\circ\gamma(t))\cdot \partial_x^a \partial_y^b
  h(\gamma_i(t))}
  {\lambda^{(1+b)n}\ell^{(1+a)n}}dt\\
  &= \sum_{1\le i\le \ell^n}\sum_{a+b\le \rho-1} \int \frac{\varphi(
  g_i^{-1}(s))\cdot Q_{a,b}(\pi\circ\gamma\circ g_i^{-1}(s))\cdot (
  g_i^{-1})'(s)\cdot \partial_x^a \partial_y^b
  h(\hat{\gamma}_i(s))}{\lambda^{(1+b)n}\ell^{(1+a)n} }  ds.
  \end{align*}
Since the $C^{a+b}$ norm of the function $ s\mapsto \varphi(
g_i^{-1}(s))\cdot Q_{a,b}(\pi\circ\gamma\circ g_i^{-1}(s))\cdot (
g_i^{-1})'(s)$ is bounded by  some constant (depending on $n$) from
(\ref{eqn:distg}), we have
  \begin{equation}\label{eqn:Psi}
  \left|\int \varphi(t) \Psi(\gamma(t)) dt \right|\le C(n)
  \|h\|^\ac_{\rho-1}
  \end{equation}
where $C(n)$ may depend on $n$ but not on $h$.

The first integral on the right hand side of (\ref{eqn:p}) is
written as
  \begin{align*}
  &\int \varphi(t) \Phi(\gamma(t)) dt =\sum_{1\le i\le \ell^n}
  \sum_{k=0}^{{\alpha}}\int  \varphi(t) \frac{ Q_k(\pi\circ
  \gamma_i(t))\cdot \partial_x^{{\alpha}-k}\partial_y^{\beta+k}
  h(\gamma_i(t))}
  {\lambda^{(1+\beta+k)n}\ell^{(1+{\alpha}-k)n}} dt\\
  &=\sum_{1\le i\le \ell^n} \sum_{k=0}^{{\alpha}}\int \frac{
  \varphi( g_i^{-1}(s))\cdot Q_k(\pi\circ \gamma\circ
  g_i^{-1}(s))\cdot ( g_i^{-1})'(s)
   \cdot \partial_x^{{\alpha}-k}\partial_y^{\beta+k}
  h(\hat{\gamma}_i(s))}
  {\lambda^{(1+\beta+k)n}\ell^{(1+{\alpha}-k)n}} ds.
  \end{align*}
For a while, we fix $1\le i\le \ell^n$. Since
  \begin{align*}
  \frac{d}{ds}&(\partial_x^{{\alpha}-k}\partial_y^{\beta+k-1}
  h(\hat{\gamma}_i(s)))\\
  &=(\pi\circ \hat{\gamma}_i)'(s) \cdot
  \partial_x^{{\alpha}-k+1}\partial_y^{\beta+k-1}
  h(\hat{\gamma}_i(s)) +\partial_x^{{\alpha}-k}\partial_y^{\beta+k}
  h(\hat{\gamma}_i(s)),
  \end{align*}
integration by part yields, for any $\psi\in
C^\rho(\mathcal{D}(\hat\gamma_i))$,
  \begin{align*}
  \int \frac{d\psi}{ds}(s)\cdot  \frac{
  \partial_x^{{\alpha}-k}\partial_y^{\beta+k-1}
  h(\hat{\gamma}_i(s))}
  {\lambda^{(1+\beta+k)n}\ell^{(1+{\alpha}-k)n}} ds &= - \int
  \tilde{\psi}(s)
  \frac{\partial_x^{{\alpha}-k+1}\partial_y^{\beta+k-1}
  h(\hat{\gamma}_i(s)) }
  {\lambda^{(1+\beta+k-1)n}\ell^{(1+{\alpha}-k+1)n}}ds\\
  &\qquad\qquad - \int \psi(s)
  \frac{\partial_x^{{\alpha}-k}\partial_y^{\beta+k}
  h(\hat{\gamma}_i(s)) }
  {\lambda^{(1+\beta+k)n}\ell^{(1+{\alpha}-k)n}}ds
  \end{align*}
where $ \tilde{\psi}(s)=\lambda^{-n}\ell^n(\pi\circ
\hat\gamma_i)'(s) \psi(s)$. This implies
  \begin{equation}\label{eqn:rec}
  \begin{aligned}
  \left|\int \psi(s) \frac{
  \partial_x^{{\alpha}-k}\partial_y^{\beta+k} h(\hat{\gamma}_i(s))}
  {\lambda^{(1+\beta+k)n}\ell^{(1+{\alpha}-k)n}} ds\right| &\le
  \left|\int \tilde{\psi}(s)
  \frac{\partial_x^{{\alpha}-k+1}\partial_y^{\beta+k-1}
  h(\hat{\gamma}_i(s)) }
  {\lambda^{(1+\beta+k-1)n}\ell^{(1+{\alpha}-k+1)n}}ds\right|\\
  &\qquad\qquad\qquad +C(n)\|\psi\|_{C^{\rho}}\|h\|^\ac_{\rho-1}
  \end{aligned}
  \end{equation}
where $C(n)$  may depend on $n$ but not on $h$ nor $\psi$. Put
  \[
  \psi_0(s)=\varphi( g_i^{-1}(s))\cdot {Q_k(\pi\circ\gamma\circ
  g_i^{-1}(s))}\cdot{( g_i^{-1})'(s)}.
  \]
By using  the last inequality repeatedly, we  obtain
  \begin{align*}
  \left|\int \frac
  {\psi_0(s)\partial_x^{\alpha-k}\partial_y^{\beta+k}
  h(\hat{\gamma}_i(s))} {\lambda^{(1+\beta+k)n}\ell^{(1+\alpha-k)n}}
  ds\right|&\le \left|\int
  \frac{\psi_{\beta+k}(s)\partial_x^{\rho}h(\hat{\gamma}_i(s))}
  {\lambda^n \ell^{(1+\rho)n}}
  ds\right|\\
  &\qquad\qquad
  +\sum_{j=0}^{\beta+k-1}C(n)\norm{\psi_{j}}_{C^{\rho}}\|h\|^\ac_{\rho-1}
  \end{align*}
where $ \psi_j(s)=\lambda^{-nj}\ell^{nj}((\pi\circ
\hat\gamma_i)'(s))^j \psi_0(s) =\lambda^{-nj} ((\pi\circ
\gamma\circ  g_i^{-1})'(s))^j \psi_0(s)$. Since
$\|\psi_j\|_{C^\rho}<C_0\lambda^n$ for $0\le j\le \beta+k$ for
some constant $C_0$ from (\ref{eqn:distg}),  we get
  \begin{align*}
  \left|\int \psi_0(s) \frac
  {\partial_x^{\alpha-k}\partial_y^{\beta+k} h(\hat{\gamma}_i(s))}
  {\lambda^{(1+\beta+k)n}\ell^{(1+\alpha-k)n}} dt\right|&\le C_0
  \ell^{-(1+\rho)n}\|h\|^\ac_\rho +C(n)\|h\|^\ac_{\rho-1}.
  \end{align*}
Summing up this inequality for $\gamma_i$, $1\le i\le \ell^n$, we
obtain
  \[
  \left|\int\varphi(t)\Phi(\gamma(t)) dt \right| \le C_0 \ell^{-\rho
  n}\|h\|^\ac_\rho +C(n)\|h\|^\ac_{\rho-1}
  \]
for some constant $C_0$. This and (\ref{eqn:Psi}) give
(\ref{eqn:GLLY1}). The proof of (\ref{eqn:GLLY2}) is obtained in a
similar but much simpler manner.
\end{proof}

\def\zero{\mathbf{0}}

\section{Main Lasota-Yorke inequality}\label{sec:ly}

In this section, we prove the following proposition.

\begin{proposition}
\label{prop:mainLY}
There exists a constant $B_0$ independent of $q$ and a constant
$C(q)$ such that, for all $\varphi\in C^{r}(\ann)$, for all integer
$\ru$ with $s+1<\ru\leq r-1$,
  \begin{equation*}
  \norm{P^q \varphi}_{W^s}^2 \leq \frac{B_0 \e(q)}{(\lambda^{1+2s}
  \ell)^q} \norm{\varphi}_{W^s}^2
  + C(q) \norm{\varphi}_{W^s} \norm{\varphi}_{\ru}^\ac.
  \end{equation*}
\end{proposition}

First of all, we introduce some notation and prove some elementary
facts concerning the Sobolev norm $\norm{\cdot}_{W^s}$. The
Fourier transform $\F\varphi$ of $\varphi\in C^{r}(\ann)$ is a
function on $\Z\times \R$ defined by
  \[
  \F\varphi(\xi,\eta)=\frac{1}{\sqrt{2\pi}}\int_{\ps}
  \varphi(x,y)\exp\left(-\ic(2\pi \xi x+\eta y)\right) dx dy.
  \]
For $s\ge 0$ and for $\varphi_1, \varphi_2\in C^{r}(\ann)$,  we
define
  \[
  (\varphi_1, \varphi_2)_{W^s}:=(\varphi_1,
  \varphi_2)^*_{W^s}+(\varphi_1, \varphi_2)_{L^2}
  \]
where
  \[
  (\varphi_1, \varphi_2)_{W^s}^*:=
  \sum_{\xi=-\infty}^{\infty}\int_{\R}\F \varphi_1 (\xi, \eta) \cdot
  \overline{\F\varphi_2 (\xi, \eta)}\cdot ((2\pi\xi)^2+\eta^2)^{s}
  d\eta.
  \]
The Sobolev norm is defined by $\|\varphi\|_{W^s}=\sqrt{(\varphi,
\varphi)_{W^s}}$. Note that we have
\begin{equation}\label{eqn:Sob1}
(\varphi_1, \varphi_2)_{W^s}^*= \sum_{\alpha+\beta= [s]}
b_{{\alpha}\beta} (\partial_x^{\alpha}\partial_y^\beta\varphi_1,
\partial_x^{\alpha}\partial_y^\beta\varphi_2)^*_{W^{s-[s]}}
\end{equation}
where  $b_{{\alpha}\beta}$ are positive integers satisfying $
(X^2+Y^2)^{[s]}=\sum_{{\alpha},\beta}b_{\alpha\beta}
X^{2\alpha}Y^{2\beta}$. Especially, if $s$ is an integer, we have
  \begin{equation}\label{eqn:Sob2}
  (\varphi_1, \varphi_2)_{W^s}^*= \sum_{\alpha+\beta=
  s}b_{{\alpha}\beta}\int_{\ps}
  \partial_x^{\alpha}\partial_y^\beta\varphi_1 (x,y) \cdot
  \overline {\partial_x^{\alpha}\partial_y^\beta\varphi_2 (x,y)}  dxdy.
  \end{equation}
In case $s$ is not an integer, we shall use the following formula
(\cite[pp 240]{H}): there exists a constant $B>0$ that depends
only on $0<\sigma<1$ such that
  \begin{align}\label{eqn:n}
  &(\varphi_1, \varphi_2)_{W^\sigma}^*=\\
  &B\int_{\ps}\!\!\!\!\!dxdy \int_{\R^2}
  \frac{(\varphi_1(x+u,y+v)-\varphi_1(x,y))\overline {(\varphi_2(x+u,y+v)-
  \varphi_2(x,y))}}{(u^2+v^2)^{1+\sigma}}dudv.
  \notag
  \end{align}
\begin{lemma}\label{lm:Sob}
{\rm(1)} For $0\le t<s\le r$ and $\epsilon>0$, there is a constant
$C(\epsilon,t,s)$ such that
  \begin{equation*}\label{eqn:ts}
  \|\varphi\|_{W^t}^2 \leq \epsilon\|\varphi\|_{W^s}^2+C(\epsilon,
  t,s)\|\varphi\|_{L^1}^2\quad \mbox{for $\varphi\in C^r(\ann)$.}
  \end{equation*}
{\rm(2)} For $\epsilon>0$, there exists a constant $C(\epsilon,s)$
with the following property: if the supports of functions
$\varphi_1, \varphi_2\in C^r(\ann)$ are disjoint and the distance
between them is greater than $\epsilon$, it holds
  \begin{equation*}\label{eqn:disj}
  \bigl| (\varphi_1, \varphi_2)_{W^s}\bigr|\le
  C(\epsilon,s)\|\varphi_1\|_{L^1}\|\varphi_2\|_{L^1}.
  \end{equation*}
\end{lemma}
\begin{proof} (1) follows from  the definition of the norm  and the
fact $\|\F\varphi\|_{L^\infty}\le \|\varphi\|_{L^1}$. If $s$ is an
integer, (2) is trivial since $(\varphi_1, \varphi_2)_{s}=0$ by
(\ref{eqn:Sob2}). Suppose that $s$ is not an integer.
Using (\ref{eqn:Sob1}) and (\ref{eqn:n}) with the assumption on
the disjointness of the supports and changing variables, we can
rewrite $(\varphi_1, \varphi_2)_{W^s}^*$ as
  \[
  -2B\sum_{\alpha+\beta=
  [s]}\int_{\ps}dxdy\int_{\R^2}\frac{b_{\alpha\beta}\cdot
  \partial_x^\alpha\partial_y^\beta\varphi_1(x+u, y+v)\cdot
  \overline {\partial_x^\alpha\partial_y^\beta\varphi_2(x,y)}}
  {(u^2+v^2)^{1+\sigma}}dudv
  \]
where $\sigma=s-[s]$. Integrating $[s]$ times by part on $(u,v)$,
then changing variables and integrating again $[s]$ times by part,
we obtain
  \[
  (\varphi_1,
  \varphi_2)_{W^s}^*=\int_{\ps}dxdy\int_{\R^2}\frac{\varphi_1(x+u,
  y+v) \overline {\varphi_2(x,y)}
\tilde{B}(u,v)}{(u^2+v^2)^{1+\sigma+2[s]}}dudv
  \]
where $\tilde{B}(u,v)$ is a polynomial of $u$ and $v$ of order
$2[s]$. With this and the assumption, we can conclude the
inequality in (2).
\end{proof}

The norm $\norm{\cdot}^\ac$ will be used through the following
lemma. Let $\C^*$ be the cone in $\R^2$ defined by
  \[
  \C^*=\{(\xi,\eta)\in \R^2\tq |\eta|\le {\alpha_0}^{-1}|\xi|\},
  \]
so that
$DT^*_{\x}(\C^*)\subset \C^*$ for  $\x\in \ps$.
\begin{lemma}\label{lm:F}
Let $\ru$ be an integer with $s+1<\ru\leq r-1$.
Let $\a$ and $\c$ elements
of $\sym^q$ and $\sym^p$ respectively, and $\chi:\ps\to \R$ a
$C^\infty$ function supported on $\P_*(\c\a)\times \R$. Take
$(\xi,\eta)\in \Z\times\R\backslash\{(0,0)\}$
such that, for any $\x\in \P_*(\c\a)\times \R$,
$(DT^q_{\x})^*(\xi,\eta)\in \C^*$.
Then, for any $\varphi \in C^r$,
\begin{equation}\label{eqn:F}
\left|(\xi^2+\eta^2)^{\ru/2}\F (P^q(\chi\cdot \varphi))(\xi,
\eta)\right|\leq C(q,\chi) \norm{\varphi}_{\ru}^\ac,
\end{equation}
where $C(q,\chi)$ may depend on $q$ and $\chi$.
\end{lemma}
\begin{proof} Let $(\xi, \eta)$ be a vector satisfying the assumption.
Let $\Gamma$ be the set of line segments on $\ps$ that are the
intersection of a line  normal to $(\xi, \eta)$ with the region
$\P_*(\c)\times\R$. We parametrize the segments in $\Gamma$
by length. Since the support of $P^q(\chi\cdot
\varphi)$ is contained in $\ann\cap (\P_*(\c)\times\R)$, the left
hand side of \eqref{eqn:F} is bounded by some constant multiple of
  \begin{equation}\label{eqn:curve}
  \sup_{\gamma\in \Gamma} \int_{\gamma}
  \partial^{\ru}P^q(\chi\cdot \varphi) dt
  \end{equation}
where $\partial$ is partial derivative with respect to $x$ if
$|\xi|>|\eta|$ and that with respect to $y$ otherwise. For each
$\gamma\in \Gamma$, there exists a unique backward image
$\tilde\gamma$ of  $T^q$ that is contained in $\P_*(\c\a) \times \R$.
If $\x\in \tilde \gamma$ and $u$ is tangent to $\gamma$ at $T^q(\x)$, then
  \begin{equation*}
  0= \langle u, (\xi,\eta) \rangle
  = \langle (DT^q_{\x})^{-1}u, (DT^q_{\x})^*(\xi,\eta) \rangle.
  \end{equation*}
By assumption, $(DT^q_{\x})^*(\xi,\eta) \in \C^*$, whence
$(DT^q_{\x})^{-1}u \in \C$ (by definition \eqref{eq:defC} of
$\C$). Hence, $\tilde\gamma$ is the composition
$\hat\gamma\circ \psi$ of an element
$\hat\gamma$ of $\Omega$ and a $C^r$ diffeomorphism $\psi$. By obvious
estimates on the distortion of $T^m$ for $0\le m \le q$
and by the definition of the norm
$\|\cdot\|^\ac_{\ru}$, we get that \eqref{eqn:curve}
 is
bounded by $C \norm{\varphi}_{\ru}^\ac$.
\end{proof}

Let $\{\chi_{\c}:S^1\to \R\}_{\c\in \sym^p}$ be a
$C^\infty$ partition of unity subordinate to the covering $\{\inte
\P_*(\c)\}_{\c\in \sym^p}$, whence $\supp(\chi_{\c})\subset \inte
\P_*(\c)$.  Define a function $\chi_{\c\a}$ by
$\chi_{\c\a}(\tau_{\c,\a}^{-q}x) =\chi_{\c}(x)$ if $x\in \P_*(\c)$,
and extend it by $0$ elsewhere. Then the functions $\chi_{\c\a}$
for $(\a,\c)\in \sym^q\times\sym^p$ are again a
$C^\infty$ partition of unity.  To keep the notation simple, we will still
use $\chi_\c$ and $\chi_{\c\a}$ to denote $\chi_\c \circ \pi$ and
$\chi_{\c\a} \circ \pi$.
\begin{lemma}\label{lm:sum}
There is a constant $C>0$ such that, for any $\varphi\in
C^r(\ann)$, it holds
\begin{equation}\label{eqn:sum1}
\sum_{(\a,\c)\in
\sym^q\times\sym^p}\|\chi_{\c\a}\varphi\|_{W^s}^2\le
2\|\varphi\|_{W^s}^2+C\|\varphi\|_{L^1}^2
\end{equation}
and
\begin{equation}\label{eqn:sum2}
\|\varphi\|_{W^s}^2\le 7\sum_{\c\in
\sym^p}\|\chi_{\c}\varphi\|_{W^s}^2 +C\|\varphi\|_{L^1}^2.
\end{equation}
\end{lemma}
\begin{proof} Since the  claims are obvious when $s=0$,  we assume $s>0$. 
Let $t$ be the largest integer that is (strictly) less than $s$.
Then for every $\epsilon>0$ we have
  \[ \sum_{(\a,\c)\in
  \sym^q\times\sym^p}\|\chi_{\c\a}\varphi\|_{W^s}^2\le
  (1+\epsilon)\|\varphi\|_{W^s}^2+C(\epsilon)\|\varphi\|_{W^{t}}^2.
  \]
Indeed, we can check this by using
(\ref{eqn:Sob2}) if $s$ is an integer and by using
(\ref{eqn:Sob1}) and (\ref{eqn:n}) instead of (\ref{eqn:Sob2})
otherwise. Hence (\ref{eqn:sum1}) follows from lemma
\ref{lm:Sob}(1).

 From lemma \ref{lm:Sob}(2), we have $(\chi_{\c}\varphi,
\chi_{\c'}\varphi)_{W^s}\le C\|\chi_{\c}\varphi\|_{L^1}\|
\chi_{\c'}\varphi\|_{L^1}\le C\|\varphi\|_{L^1}^2$ for some
constant $C>0$ if the closures of $\P_*(\c)$ and $\P_*(\c')$
do not intersect. Also we have $(\chi_{\c}\varphi,
\chi_{\c'}\varphi)_{W^s}\le (\|\chi_{\c}\varphi\|_{W^s}^2
+\|\chi_{\c'}\varphi\|_{W^s}^2)/2$ in general.
Applying these to $\|\varphi\|_{W^s}^2=\sum_{(\c,\c')\in
\sym^p\times\sym^p}(\chi_{\c}\varphi, \chi_{\c'}\varphi)_{W^s}$,
we obtain (\ref{eqn:sum2}).
\end{proof}

We start the proof of Proposition \ref{prop:mainLY}.
>From \eqref{eqn:sum2}, we
have
  \begin{align*}
  \|P^q(\varphi)\|_{W^s}^2&\le 7\sum_{\c\in
  \sym^p}\|\chi_{\c}P^q(\varphi)\|_{W^s}^2+C\norm{\varphi}_{L^1}^2
  \\&
  \le 7\sum_{\c\in
  \sym^p}\left\|\sum_{\a\in \sym^q}
  P^q(\chi_{\c\a}\varphi)\right\|_{W^s}^2+C\norm{\varphi}_{L^1}^2.
  \end{align*}
So we will estimate
  \[
  \left\|\sum_{\a\in \sym^q} P^q(\chi_{\c\a}\varphi)\right\|_{W^s}^2=
  \sum_{(\a, \b)\in \sym^q\times \sym^q}(P^q(\chi_{\c\a}\varphi),
  P^q(\chi_{\c\b}\varphi))_{W^s}
  \]
for $\c\in \sym^p$.

Consider first a pair $(\a, \b)\in \sym^q\times
\sym^q$ such that $\a\pitchfork_{\c}\b$. For any $(\xi,\eta)\in
\Z\times \R\backslash \{(0,0)\}$,
this implies that either $(DT_{\x}^q)^*(\xi,\eta)\in \C^*$ for all
$\x\in \P_*(\c\a)\times \R$, or $(DT_{\x}^q)^*(\xi,\eta)\in \C^*$ for all
$\x\in \P_*(\c\b)\times \R$. Let $U$ be the set of all $(\xi,\eta)\in \Z\times\R$ such
that the first possibility holds, and $V=(\Z\times\R) \backslash U$. If
$(\xi,\eta)\in U$,
by Lemma \ref{lm:F}, there exists a constant $C>0$
such that $\left|\F
P^q(\chi_{\c\a}\varphi)(\xi,\eta)\right| \leq C (\xi^2+\eta^2)^{-\ru/2}
\norm{\varphi}_{\ru}^\ac$. Moreover, $\left| \F
P^q(\chi_{\c\a}\varphi)(\xi,\eta) \right| \leq C
\norm{\varphi}_{L^1}$, which is bounded by $C\norm{\varphi}_{\ru}^\ac$
by \eqref{ac_greater_L1}.
Hence, $\left|\F
P^q(\chi_{\c\a}\varphi)(\xi,\eta)\right| \leq C
(1+\xi^2+\eta^2)^{-\ru/2} \norm{\varphi}_{\ru}^\ac$.
So we have, for some constant $C$,
  \begin{align*}
  &\left|\; \sum_{\xi=-\infty}^{\infty} \int \mathbf{1}_U(\xi,\eta)\cdot  (1+\xi^2+\eta^2)^s \F P^q (\chi_{\c\a} \varphi) \cdot
  \overline{ \F P^q (\chi_{\c\b} \varphi)}\; d\eta\;\right|
  \\
  &\quad 
  \leq C \left(\sum_{\xi=-\infty}^{\infty} \int \mathbf{1}_U(\xi,\eta)\cdot (1+\xi^2+\eta^2)^s | \F P^q (\chi_{\c\a}
  \varphi) |^2 d\eta \right)^{1/2} \norm{ P^q(\chi_{\c\b} \varphi)}_{W^s}
  \\&\quad 
  \leq C  \norm{\varphi}_{\ru}^\ac
  \norm{\varphi}_{W^s},
  \end{align*}
since the function $(1+\xi^2+\eta^2)^{-\ru+s}$ is integrable by the
assumption $s<\ru-1$. The same inequality holds on $V$, and we obtain
  \begin{equation}\label{eqn:c1}
  \left|(P^q(\chi_{\c\a}\varphi),
  P^q(\chi_{\c\b}\varphi))_{W^s}\right| \leq
  C \norm{\varphi}_{\ru}^\ac \norm{\varphi}_{W^s}.
  \end{equation}
For the sum over $\a$ and $\b$ such that $\a\not\pitchfork_\c \b$,
we have
  \begin{align}
  \sum_{\a\not\pitchfork_\c \b}(P^q(\chi_{\c\a}\varphi),
  P^q(\chi_{\c\b}\varphi))_{W^s}
  &\le \sum_{\a\not\pitchfork_\c \b}
  \frac{\|P^q(\chi_{\c\a}\varphi)\|^2_{W^s}+\|P^q(\chi_{\c\b}
  \varphi)\|^2_{W^s}}{2}\notag\\
  &\le \e(q)  \sum_{\a\in \sym^q}
  \|P^q(\chi_{\c\a}\varphi)\|^2_{W^s}.\label{eqn:c2}
  \end{align}
For the terms in the last sum, we have the estimate
  \begin{equation}\label{eqn:c3}
  \|P^q(\chi_{\c\a}\varphi)\|^2_{W^s}\le
  \frac{C_0\|\chi_{\c\a}\varphi\|_{W^s}^2}{\lambda^{(1+2s)q}\ell^{q}}
  +C\norm{\varphi}_{L^1}^2
  \end{equation}
where $C_0$ is a constant that depends only on $\lambda$, $\ell$
and $\kappa$. Indeed, we can check this by using \eqref{eqn:Sob2}
and \eqref{eqn:Tloc} if $s$ is an integer and by using
\eqref{eqn:Sob1} and \eqref{eqn:n} instead of \eqref{eqn:Sob2}
otherwise.

>From \eqref{eqn:c1}, \eqref{eqn:c2}, \eqref{eqn:c3},
\eqref{eqn:sum1} and \eqref{ac_greater_L1}, we obtain
  \begin{align*}
  \sum_{\c\in \sym^p} \left\|\sum_{\a\in \sym^q}
  P^q(\chi_{\c\a}\varphi)\right\|_{W^s}^2
  &
  \le
  \frac{C_0\e(q)}{\lambda^{(1+2s)q}\ell^q}
  \sum_{(\a,\c)\in \sym^q\times \sym^p}
  \|\chi_{\c\a}\varphi\|_{W^s}^2
  + C \norm{\varphi}_{W^s} \norm{\varphi}_{\ru}^\ac
  \\& \le
  \frac{2C_0\e(q)}{\lambda^{(1+2s)q}\ell^q}
  \norm{\varphi}_{W^s}^2
  + C \norm{\varphi}_{W^s} \norm{\varphi}_{\ru}^\ac,
  \end{align*}
and hence Proposition \ref{prop:mainLY}.

\section{Proof of the main theorems}
\label{sec:proofs_thms}

We will use Lemma \ref{lm:P} and Proposition \ref{prop:mainLY} to
study the properties of $P$ acting on the space $C^r(D)$ equipped with the norms
$\norm{\cdot}_{\ru}^\ac$ and $\norm{\cdot}_{W^s}$.

\begin{lemma}
\label{lem:delta1}
Let $\delta\in (\ell^{-1},1)$.
There exists $C>0$ such that, for integer $1\leq \rho \leq r-1$, for $n\in
\N$,
  \begin{equation*}
  \norm{P^n h}_\rho^\ac \leq C \delta^{\rho n} \norm{h}_\rho^\ac + C
  \norm{h}_{\rho-1}^\ac.
  \end{equation*}
\end{lemma}
\begin{proof}
We prove it by induction on $\rho$.
Let $\rho\geq 1$. By Lemma  \ref{lm:P},
there exists $N\in \N$ and $C>0$
such that
  \begin{equation}
  \label{eq:iter}
  \norm{P^N h}_{\rho}^\ac \leq \delta^{\rho N} \norm{h}_{\rho}^\ac+ C
  \norm{h}_{\rho-1}^\ac.
  \end{equation}
By the inductive assumption (and Lemma \ref{lm:P} in the $\rho=1$
case), $\norm{P^n h}_{\rho-1}^\ac \leq C \norm{h}_{\rho-1}^\ac$. Hence,
iterating \eqref{eq:iter} gives the conclusion.
\end{proof}

\begin{lemma}
\label{lem:delta2}
Let $\delta\in (\ell^{-1},1)$, and let $0 \leq \rv<\ru \leq r-1$ be
integers. Let $\nu(\ru,\rv)=\sum_{j=\rv+1}^{\ru}
\frac{1}{j}$. There exists $C>0$ such that, for $n\in \N$,
  \begin{equation*}
  \norm{P^n h}_{\ru}^\ac \leq C \delta^{n/\nu(\ru,\rv)} \norm{h}_{\ru}^\ac +
  C \norm{h}_{\rv}^\ac.
  \end{equation*}
\end{lemma}
\begin{proof}
Let $n$ be a multiple of $(r-1)!$.  Then
Lemma \ref{lem:delta1} implies by induction over
$\rv+1 \leq \rho\leq \ru$ that
  \begin{equation*}
  \norm{P^{(\frac{1}{\rho}+\dots+\frac {1} {\rv+1})n} h}_{\rho}^\ac \leq
  C \delta^n \norm{h}_{\rho}^\ac +C \norm{h}_{\rv}^\ac.
  \end{equation*}
For $\rho=\ru$, we obtain
$\norm{P^{\nu(\ru,\rv) n}h}_{\ru}^\ac \leq C \delta^n \norm{h}_{\ru}^\ac + C
\norm{h}_{\rv}^\ac$.
\end{proof}

\begin{theorem}
\label{Theorem:LY}
Assume that $\frac{B_0 \e(q)}{(\lambda^{1+2s}\ell)^q} <1$.  Let $0
\leq \rv<\ru \leq r-1$ be integers with $s<\ru-1$, and let
$\nu=\nu(\ru,\rv)$ be
as in the previous lemma. Let
  \begin{equation*}
  \gamma \in \left( \max\left( \ell^{-1/\nu}, \sqrt{\frac{(B_0
  \e(q))^{1/q}}{\lambda^{1+2s}\ell}} \right), 1 \right).
  \end{equation*}
Let $\norm{\varphi}:=
\norm{\varphi}_{W^s}+\norm{\varphi}_{\ru}^\ac$. There exists a
constant $C$ such that, for all $n\in \N$,
  \begin{equation*}
  \norm{P^n \varphi} \leq C \gamma^n \norm{\varphi}
  + C \norm{\varphi}_{\rv}^\ac.
  \end{equation*}
\end{theorem}
\begin{proof}
Since $\sqrt{a+b}\leq \sqrt{a}+\sqrt{b}$ and $\sqrt{ab} \leq \epsilon
a +\epsilon^{-1} b$, Proposition \ref{prop:mainLY} implies
  \begin{equation*}
  \norm{P^q \varphi}_{W^s} \leq  \left(\frac{(B_0
  \e(q))^{1/q}}{\lambda^{1+2s}\ell}\right)^{q/2}\norm{\varphi}_{W^s} +
  \epsilon \norm{\varphi}_{W^s}+
  C(\epsilon) \norm{\varphi}_{\ru}^\ac.
  \end{equation*}
Since $\left(\frac{(B_0
\e(q))^{1/q}}{\lambda^{1+2s}\ell}\right)^{q/2}< \gamma^q$, taking
$\epsilon$ small enough yields
  \begin{equation*}
  \norm{P^q \varphi}_{W^s} \leq \gamma^q \norm{\varphi}_{W^s}+C
  \norm{\varphi}_{\ru}^\ac.
  \end{equation*}
Iterating
this equation $K$ times gives
  \begin{equation}
  \label{eq:sum1}
  \norm{P^{Kq} \varphi}_{W^s} \leq \gamma^{Kq} \norm{\varphi}_{W^s} +
  C(K) \norm{\varphi}_{\ru}^\ac,
  \end{equation}
for some constant $C(K)$.
If $K$ is large enough, the choice of $\gamma$ and Lemma
\ref{lem:delta2} also yield
  \begin{equation}
  \label{eq:sum2}
  \norm{P^{Kq} \varphi}_{\ru}^\ac \leq \frac{\gamma^{Kq}}{2}
  \norm{\varphi}_{\ru}^\ac + C'(K) \norm{\varphi}_{\rv}^\ac.
  \end{equation}
Fix such a $K$, and define a norm $\norm{\varphi}^*:=
\norm{\varphi}_{W^s} + 2 C(K)\gamma^{-Kq}
\norm{\varphi}_{\ru}^\ac$. Adding
\eqref{eq:sum1} and
\eqref{eq:sum2} gives
  \begin{equation*}
  \norm{P^{Kq} \varphi}^* \leq \gamma^{Kq} \norm{\varphi}^* + C
  \norm{\varphi}_{\rv}^\ac.
  \end{equation*}
Iterating this equation (and remembering $\norm{P^n
\varphi}_{\rv}^\ac \leq C \norm{\varphi}_{\rv}^\ac$ for some
constant $C$ independent of $n$, by Lemma \ref{lem:delta1}), we
obtain the conclusion of the theorem for the norm
$\norm{\cdot}^*$. Since it is equivalent to the original norm
$\norm{\cdot}$, this concludes the proof.
\end{proof}

\begin{corollary}
\label{cor:main1}
If $B_0 \e(q) < (\lambda^{1+2s}\ell)^q$, the conclusion of Theorem
\ref{th:main1} holds for the transformation $T$.
\end{corollary}
\begin{proof}
Take $\ru=r-1$ and $\rv=0$. They satisfy the assumptions of
Theorem~\ref{Theorem:LY} since $s<r-2$.

We fix a non-negative function $\Psi_0\in
C^r(\ann)$ such that $\int \Psi_0 d\vol=1$. Put $\nu_0=\Psi_0\cdot
\vol$ and $\Psi_n=P^n\Psi_0$ for $n\ge 1$. From \eqref{eqn:ch},
the density of $T^n_*\nu_0$ is $\Psi_n$. As the sequence
$T^n_*\nu_0$ converges to the SBR measure $\mu$ for $T$ weakly, we
have
  \begin{equation}\label{eqn:wc}
  \lim_{n\to \infty}\left(\Psi_n,\varphi\right)_{L^2}=\int \varphi
  d\mu
  \end{equation}
for any continuous function $\varphi$ on $\ps$ with compact
support. By Theorem~\ref{Theorem:LY}, the sequence $\Psi_n$ for
$n\ge 1$ is bounded with respect to the norm $\norm{\cdot}$, hence
also for the norm $\|\cdot \|_{W^s}$. Then there is a subsequence
$n(i)\to \infty$ such that $\Psi_{n(i)}$ converges weakly to some
element $\Psi_\infty$ in the Hilbert space $W^s(\ps)$. This and
\eqref{eqn:wc} imply $ \int \Psi_\infty \varphi d\vol=\int \varphi
d\mu$ for any continuous function $\varphi$ on $\ps$ with compact
support. Thereby the density of the SBR measure $\mu$ is
$\Psi_\infty\in W^s(\ps)$.
\end{proof}

\begin{corollary}
\label{cor:main2}
Let $1/2<s<r-2$.
Assume that $\frac{B_0 \e(q)}{(\lambda^{1+2s}\ell)^q} <1$. If
  \begin{equation*}
  \gamma \in \left( \sqrt{\frac{(B_0
  \e(q))^{1/q}}{\lambda^{1+2s}\ell}}, 1 \right),
  \end{equation*}
the conclusion of Theorem~\ref{th:main2} holds for the transformation
$T$ and this $\gamma$.
\end{corollary}
\begin{proof}
Let $\ru$ be the smallest integer such that $s<\ru-1$, and $\rv$ the
largest integer such that $\rv<s-1/2$. They satisfy the assumptions of
Theorem~\ref{Theorem:LY}. Moreover, $\nu(\ru,\rv)\leq
1+\frac{1}{2}+\frac{1}{3}<2$. Hence, $\ell^{-1/\nu} < \frac{1}{\sqrt{\ell}}<
\sqrt{\frac{(B_0 \e(q))^{1/q}}{\lambda^{1+2s}\ell}}$.

Let $\boB$ be the completion of $C^r(\ann)$ with respect to the
norm $\norm{\cdot}$. It is a Banach space included in $W^s(D)$ and
containing $C^{r-1}(D)$.
Theorem~\ref{Theorem:LY} gives a Lasota-Yorke inequality between
$\boB$ and the space $\boB'$ obtained by completing $C^r(\ann)$
for the norm $\norm{\cdot}_{\rv}^\ac$. Hence, the result is a
standard consequence of Hennion's Theorem \cite{He}, if we can
prove that the unit ball of $\boB$ is relatively compact in
$\boB'$.

The embedding of $\boB$ in $W^s(\ann)$ is continuous.
Let $t\in (\rv+1/2,s)$. The embedding of $W^s(\ann)$ in $W^t(\ann)$ is
compact by Sobolev's embedding theorem. To conclude, it is sufficient
to check that the injection $W^t(\ann) \to \boB'$ is
continuous. Since $t>\rv+1/2$, \cite[Theorem 7.58 (iii)]{A} (applied with
$p=q=2$ , $k=1$ and $n=2$) proves that, for any smooth curve
$\boC \subset \ann$, for any $\varphi\in W^t(\ann)$,
  \begin{equation*}
  \norm{\partial^\alpha_x \partial^\beta_y \varphi}_{L^2(\boC)}
\leq C(\boC) \norm{\varphi}_{W^t(\ann)}
  \end{equation*}
whenever $\alpha$ and $\beta$ are non-negative integers satisfying
$\alpha+\beta \leq \rv$. The constant $C(\boC)$ can be chosen
uniformly over all curves of $\Omega$, and we obtain
$\norm{\varphi}_{\rv}^\ac \leq C \norm{\varphi}_{W^t(\ann)}$.
\end{proof}

For $\beta>0, \kappa>0$ and $\lambda\in (0,1)$, let
  \[
  \mathcal{E}(\beta,\kappa, \lambda)=\left\{f \in \mathcal{U}_\kappa \;;\;
  \limsup_{q\to\infty} \frac 1 q \log \e(q)> \beta\right\}.
  \]
Note that this definition depends on $\kappa$ and $\lambda$ through
$\e(q)$, since $\e(q)$ is defined in terms of $\alpha_0=\kappa/(1-\lambda)$.

Since the quantity $\e(q)$ depends on $f\in \mathcal{U}_\kappa$ upper
semi-continuously and since we can take arbitrarily large $\kappa$
in the beginning, Theorems~\ref{th:main1}, \ref{th:main2}
and~\ref{th:main} follow
from Corollaries~\ref{cor:main1} and~\ref{cor:main2} and the next
proposition.
\begin{proposition}\label{pp:transv}
For any $\beta>0$ and $\lambda>0$, there is a finite collection of $C^{\infty}$
functions $\varphi_{i}:S^{1}\to \R$, $i=1,2,\cdots, m$ and a constant
$D_0>0$ such that,
for any $\kappa>D_0$ and any $C^{r}$ function $g\in
\mathcal{U}_{\kappa-D_0}$, the subset
  \[
  \left\{(t_{1},t_{2},\cdots,t_{m})\in [-1,1]^m
  \;\left|\;g+\sum_{i=1}^{m}
  t_{i}\varphi_{i} \in \mathcal E(\beta,\kappa, \lambda)\right.\right\}
  \]
is a Lebesgue null subset on $[-1,1]^m$.
\end{proposition}

This proposition has essentially been  proved in \cite{T}.
For completeness, we give a proof of it in the next
section.

\section{Genericity of the transversality condition}\label{sec:gene}

In this section, we give a proof of Proposition~\ref{pp:transv}.
For  a $C^{2}$ function $g$  and  $C^{\infty}$ functions
$\varphi_{i}$, $1\le i\le m$, on $S^{1}$, we consider a family of
functions
\begin{equation}\label{eqn:fam}
f_{\t}(x)=g(x)+\sum_{i=1}^{m}t_{i}\varphi_{i}(x): S^{1}\to \R
\end{equation}
and the  corresponding family of maps
\begin{equation}\label{eqn:Tfam} T_{\t}:\ps\to\ps,\qquad
 T_{\t}(x,y)=(\lap x, \lambda y+f_{\t}(x))
\end{equation}
with parameters $\t=(t_{1},t_{2},\cdots,t_{m})\in [-1,1]^m\subset
\R^{m}$. Put
\begin{equation} \label{eqn:S2}
S(x,\a;\t)=\sum_{i=1}^{n}\lambda^{i-1}f_{\t}([\a]_{i}(x))
\end{equation}
for $\t\in[-1,1]^m$ and a word $\a\in\sym^{n}$ of length $1\le
n\le \infty$. For a point $x\in S^{1}$ and a sequence
$\sigma=(\a_{0},\a_{1},\cdots,\a_{k})$ of elements in
$\sym^{\infty}$, we consider an affine map $
G_{x,\sigma}:\R^{m}\to\R^{k}$ defined by
\begin{equation}\label{gsigma}
G_{x,\sigma}(\t)=\left(\frac{d}{d x}S(x,\a_{i};\t)-\frac{d}{d
x}S(x,\a_{0};\t)\right)_{i=1,2,\cdots,k}.
\end{equation}
If the affine map $G_{x,\sigma}$ is surjective, we define its
Jacobian  by
\[
\Jac(G_{x,\sigma})=\frac{\vol_{k}([0,1]^k)}{\vol_{k}(G_{x,\sigma}^{-1}([0,1]^k)\cap
\Ker(G_{x,\sigma})^{\perp} )}
\]
where $\vol_k$ is the $k$-dimensional Hausdorff measure and
$\Ker(G_{x,\sigma})^{\perp}$ is the orthogonal complement of the
kernel of the linear part of $G_{x,\sigma}$, whence
\begin{equation}\label{eqn:Jac}
\vol(G_{x,\sigma}^{-1}(Y)\cap [-1,1]^m)\le
C_0\frac{\vol(Y)}{\Jac(L)}\quad \mbox{for any Borel subset
$Y\subset \R^k$}
\end{equation}
where  $C_0$ is a constant that depends only on the dimensions $m$
and $k$.

For $0<\gamma\le 1$, $\delta>0$ and $n\ge 1$, we say that the
family $T_{\t}^n$ is {\em $(\gamma,\delta)$-generic} if the
following property holds: for any finite sequence
$\{\a_{i}\}_{i=0}^d$ in $\sym^{\infty}$ such that $[\a_i]_n$ are
mutually distinct, for any $x\in S^{1}$ and for any integer $0<k<
\gamma d$, we can choose a subsequence
$\sigma=(\b_{0},\b_{1},\cdots,\b_{k})$ of length $k$ among
$\{\a_{i}\}_{i=0}^d$ so that $G_{x,\sigma}$ is surjective and
satisfies $\Jac(G_{x,\sigma})>\delta$. It is proved  in \cite{T}
that

\begin{proposition}[\cite{T}, {Proposition 15}]\label{prop:para}
For given $0<\lambda<1$, $\ell\ge 2$ and $n\ge 1$, there exists a
finite collection of $C^{\infty}$ functions $\varphi_{i}$, $1\le
i\le m$, such that  the corresponding family  $T_{\t}^{n}$ is
$(1/(n+1),1/2)$-generic, regardless of the $C^{2}$ function $g$.
\end{proposition}

Recall that we are considering fixed $\lambda$ and $\ell$. Let
$\beta>0$ be the  positive number in the  statement of Proposition
\ref{pp:transv}. We  can and do take integers $N_0\ge 2$, $d_0\ge
2$ and $n_0\ge 1$  such that
\begin{equation}\label{eqn:condq}
\lambda^{N_0-1}\ell^2<1, \quad d_0/(n_0+1)>N_0+1 \quad \mbox{and
}\quad (d_0+1)\exp(-\beta n_0/2)<1/2.
\end{equation}
Let $\varphi_i$, $1\le i\le m$, be the $C^\infty$ functions in  the
conclusion of Proposition~\ref{prop:para} for these $\lambda$,
$\ell$ and $n=n_0$. Let $D_0=\sum_{i=1}^m \norm{\varphi_i}_{C^r}$. Hence,
if $g\in \mathcal{U}_{\kappa-D_0}$ and $(t_1,\dots,t_m)\in [-1,1]^m$,
then $g+\sum t_i \varphi_i \in \mathcal{U}_\kappa$.
In order to prove the conclusion of
Proposition~\ref{pp:transv}, we pick arbitrary $g\in
\mathcal{U}_{\kappa- D_0}$
and consider the family  $T_{\t}$ defined by \eqref{eqn:fam} and
\eqref{eqn:Tfam}.

For an integer $q$, we put $p(q)=[q\log(\ell/\lambda)/\log
\ell]+1$. For a word $\c$ of finite length, let $x_{\c}$ be the
left end of $\P(\c)$. We fix a word $\a_\infty\in \sym^\infty$
arbitrarily and, for any word $\a$ of finite length, we put $\bar
\a=\a \a_\infty$.
\begin{lemma}\label{lm:E}If $f_{\t}\in \mathcal{E}(\beta,\kappa)$, we
can take arbitrarily large integer $q$ such that there exist  $1+d_0$
words $\a_i$, $0\le i\le d_0$, in $\sym^q$ and a word $\c\in
\sym^{p(q)}$ satisfying
\begin{itemize}
\item[{\rm(E1)}] $|\frac{d}{d x}S(x_{\c},\bar\a_i;\t)-\frac{d}{d x}
S(x_{\c},\bar\a_j;\t)|\le 8\lambda^q\ell^{-q}{\alpha_0}$\;\;
for any $1\le i,j\le d_0$, and
\item[{\rm(E2)}] $[\a_i]_{n_0}\neq [\a_j]_{n_0}$ if $i\neq j$.
\end{itemize}
\end{lemma}
\begin{proof}
\def\tm{\tilde{q}}
By assumption, we can take an arbitrarily large $\tm$ such that
there exist a point $x\in S^1$ and subset $E\subset \sym^{\tm}$
such that
 $\# E\ge \exp(\beta \tm)$ and
\begin{equation}\label{eqn:E1}
\left|\frac{d}{d x}S(x,\a;\t)-\frac{d}{d x}S(x,\b;\t)\right|\le
4\lambda^{\tm}\ell^{-\tm}{\alpha_0}\quad\mbox{ for $\a$ and $\b$
in $E$.}
\end{equation}
For each $0\le j\le [\tm/n_0]$, we introduce an equivalence
relation $\sim_j$ on $E$ such that $\a\sim_j \b$ if and only if
$[\a]_{jn_0}=[\b]_{jn_0}$, and let 
\[
\nu(j)=\max_{\a\in E}\#\{\b\in
E\mid \b\sim_j \a\}.
\]
 Since $\nu(0)= \# E\ge \exp(\beta \tm)$
while $\nu(j)\le \ell^{\tm-jn_0}$ obviously, there exists $0\le
j\le [\tm/n_0]$ such that $\nu(j+1)<\exp(-\beta n_0/2) \nu(j)$.
Let $j_*$ be the minimum of such integers $j$ and put
$q=\tm-n_0j_*$. Then we have $\nu(j_*)\ge \exp(\beta q)$ and $q\ge
\beta \tm/(2\log\ell)$. The equivalence class $H$
w.r.t.~$\sim_{j_*}$ of maximum cardinality contains at least
$(d_0+1)$ non-empty equivalence classes w.r.t.~$\sim_{j_*+1}$,
because
\[
\nu(j_*)-(d_0+1)\nu(j_*+1)> \nu(j_*)-(d_0+1)\exp(-\beta
n_0/2)\nu(j_*)>0
\]
by (\ref{eqn:condq}). So we can take $\b\in \sym^{\tm-q}$ and
$\a_i\in \sym^q$, $0\le i\le d$, such that $\b\a_i\in H$ for $0\le
i\le d_0$ and that (E2) holds. Put $x'=\b(x)$. It follows from
(\ref{eqn:E1}) that
\begin{equation}\label{eqn:Sdiff}
\left|\frac{d}{d x}S(x',\a_i;\t)-\frac{d}{d
x}S(x',\a_j;\t)\right|\le 4\lambda^q\ell^{-q}{\alpha_0}\quad\mbox{
for $0\le i,j\le d_0$}.
\end{equation}
Take $\c\in \sym^{p(q)}$ such that $x'\in \P(\c)$. Since the
distance between $x_{\c}$ and $x'$ is bounded by $\ell^{-p(q)}\le
\lambda^q/\ell^q$, the condition (E1) follows from
(\ref{eqn:Sdiff}) and (\ref{eqn:alpha}).
 \end{proof}
Let $\mathcal{B}^q$ be the set of pairs $(\sigma, \c)$ of a
sequence $\sigma=(\b_i)_{i=0}^{N_0}$ in $\sym^q$ and $\c\in
\sym^{p(q)}$ such that $\Jac(G_{x_\c,\bar \sigma})>1/2$, where
$\bar\sigma=({\bar{\b}}_i)_{i=0}^{N_0}$.
For $(\sigma,\c)\in \mathcal{B}^q$ with $\sigma=(\b_i)_{i=0}^{N_0}$,
we put
  \[
  Y(\sigma,\c)=G_{x_\c,\bar \sigma}^{-1}([-8(\lambda/\ell)^q
  \alpha_0,8(\lambda/\ell)^q\alpha_0]^{N_0})
  \]
and $Y(q):=\bigcup_{(\sigma,\c)\in \mathcal{B}^q} Y(\sigma,\c)$.
Since the family $T^{n_0}_{\t}$ is $(1/(n_0+1), 1/2)$-generic, the
conclusion of Lemma~\ref{lm:E} and the second condition in
\eqref{eqn:condq} imply that, if $f_{\t}\in \mathcal{E}
(\beta,\kappa,\lambda)$,
the parameter $\t$ is contained in $Y(q)$ for infinitely many $q$.
Using \eqref{eqn:Jac} and the simple estimate $\#\mathcal{B}^q\le
\ell^{q (N_0+1)+p(q)}$, we get
  \[
  \vol\left(Y(q)\right)\le C\ell^{q (N_0+1)+p(q)}(\lambda/\ell)^{qN_0}
  \]
for some constant $C>0$ . By the first condition in
\eqref{eqn:condq}, the left hand side converges to $0$
exponentially fast as $q\to \infty$. Therefore we obtain the
conclusion of Proposition~\ref{pp:transv} by Borel-Cantelli
lemma.

\end{document}